\documentclass[12pt]{article}
\def\RA{\Rightarrow}
\usepackage{latexsym}

\begin{document}
\thispagestyle{empty}

\vskip 20pt
\begin{center}
{\Large On Generalized Van der Waerden Triples}
\vskip 20pt
Bruce Landman\\
Department of Mathematical Sciences\\
University of North Carolina at Greensboro\\
Greensboro, NC 27402\\
email:  {\tt bmlandma@uncg.edu}
\vskip 10pt
and
\vskip 10pt
Aaron Robertson\\
Department of Mathematics\\
Colgate University\\
Hamilton, NY 11346\\
email:  {\tt aaron@math.colgate.edu}
\end{center}
\vskip 50pt
\begin{abstract}
Van der Waerden's classical theorem on arithmetic progressions states that
for any positive integers $k$ and $r$, there exists
a least positive integer, $w(k,r)$, such that
any $r$-coloring of $\{1,2,\dots,w(k,r)\}$ must
contain a monochromatic $k$-term arithmetic
progression $\{x,x+d,x+2d,\dots,x+(k-1)d\}$.
We investigate the following generalization of
$w(3,r)$.
For fixed positive integers $a$ and $b$ with $a \leq b$, define
$N(a,b;r)$ to be the least positive integer, if it
exists, such that any $r$-coloring of
$\{1,2,\dots,N(a,b;r)\}$ must contain a monochromatic
set of the form $\{x,ax+d,bx+2d\}$.  We show that
$N(a,b;2)$ exists if and only if $b \neq 2a$, and provide upper
  and lower bounds for it.
We then show that for  a large class of pairs
$(a,b)$,  $N(a,b;r)$ does not exist for $r$ sufficiently large. We also
give a  result on sets of the form
$\{x,ax+d,ax+2d,\dots,ax+(k-1)d\}$.
\end{abstract}
\vskip 50pt
\begin{center}
{\bf 1. Introduction}
\end{center}
B.L. van der Waerden [6] proved that for any positive 
integers $k$ and $r$, there exists a least positive
integer, $w(k,r)$, such that any $r$-coloring of
$[1,w(k,r)]=\{1,2,\dots,w(k,r)\}$ must contain a
monochromatic $k$-term arithmetic progression
$\{x,x+d,x+2d,\dots,x+(k-1)d\}$.  The only known non-trivial values of
$w(k,r)$ are $w(3,2)=9$, $w(4,2)=35$, $w(5,2)=178$, $w(3,3)=27$ and
$w(3,4)=76$. The function $w(k,r)$ is sometimes called the {\em Ramsey function}
for the collection of arithmetic progressions. In [1] the authors considered
a generalization of van der Waerden's theorem, by considering, for a given
function $f:{\bf N} \rightarrow {\bf N}$, the Ramsey function corresponding
to the collection of arithmetic progressions $\{a,a+d,a+2d,...,a+(k-1)d\}$
with the property that $d \geq f(a)$.  The Ramsey functions for other
``substitutes'' for the set of arithmetic progressions were studied in
[2], [4], and [5]. In this paper we consider a new generalization
of $w(k,r)$. To help describe this generalization, we begin with three
definitions.

\vspace{.07in}

\noindent
{\bf Definition 1.1:}  Fix $1 \leq a \leq b$.  
A set, $S$, of three natural numbers
is called an $(a,b)${\em -triple} if there exist
natural numbers $x$ and $d$ such that $S=\{x,ax+d,bx+2d\}$.

\vspace{.07in}

\noindent
{\bf Definition 1.2:}  Fix $1 \leq a \leq b$.  
Define $N(a,b;r)$ to be the least positive integer, 
if it exists, such that any $r$-coloring of 
$[1,N(a,b;r)]$ must contain a monochromatic $(a,b)$-triple.

\vspace{.07in}

\noindent
{\bf Definition 1.3:}  Fix $1 \leq a \leq b$. Define $(a,b)$ to be
{\em regular} if $N(a,b;r)$ exists for all positive integers $r$. If
$(a,b)$ is not regular, 
 the {\it degree of regularity of $(a,b)$} is
the largest $r$ such that $N(a,b;r)$ exists.
 Denote this by $dor(a,b)$.

 \vspace{.07in}

We note here that $N(1,1;r)$ is the van der Waerden 
number $w(3,r)$ so that $N(a,b;r)$ is a
generalization of $w(3,r)$, and obviously $(1,1)$ is regular.
We now discuss the
sections which follow.

In Section 2 we consider $r=2$. We
 show that, except for the case in which $b=2a$, $N(a,b;2)$
does exist;  we also 
find upper and lower bounds on $N(a,b;2)$ (for $b \neq 2a$).
  For certain pairs $(a,b)$, we obtain stronger bounds; in particular,
we use a result of [1] to deal with $N(a,2a-1;2)$ (when $a=1$ this is
just $w(3,2)$). In Section 3 we establish that $(a,b)$ is not regular for
 a rather large class
of pairs $(a,b)$, and give an upper bound (for these pairs) on
 the degree of regularity.
  We then give lower bounds 
on $N(a,b;r)$
for all $1 \leq a \leq b$ and $r>2$.
In Section 4  we make
some observations about monochromatic sets 
of the form
$\{x,ax+d,ax+2d,\dots,ax+(k-1)d\}$ for $a \geq 1$.
We establish that for $a>1$ and $k$ sufficiently large
(dependent upon $a$), we can $4$-color the
natural numbers so that no monochromatic
such $k$-set exists (this is in contrast to  van der Waerden's
theorem which says that there are arbitrarily long monochromatic 
arithmetic
progressions in any $r$-coloring of the natural
numbers). 
\vskip 50pt
\begin{center}
{\bf 2. Using Two Colors}
\end{center}
Our first theorem 
categorizes those $(a,b)$ pairs for which $N(a,b;2)$ exists, i.e., those
pairs for which $dor(a,b) \geq 2$. It also provides an upper bound on
$N(a,b;2)$ whenever it exists.

\vspace{.1in}

\noindent
{\bf Theorem 2.1:}  Let $a,b \in {\bf N}$ with $a \leq b$. Then
$dor(a,b)=1$ if and only if $b=2a$. 
 
\noindent Furthermore, if $b \neq 2a$,

\centerline{$
N(a,b;2) \leq
\left\{
\begin{array}{ll}
4a(b^3+b^2-3b-3)+2b^3+6b^2+6b&\mathrm{for}\,\, b > 2a\\
4a(b^3+2b^2+2b)-4b^2&\mathrm{for}\,\, b < 2a\\
\end{array}
\right.
$}
\vspace{.07in}

\noindent
{\em Proof.}
We first consider the case in which $b=2a$. To show that $N(a,2a;2)$ does
not exist, we exhibit a 2-coloring of {\bf N} which avoids monochromatic
$(a,2a)$-triples. Namely, color the natural numbers so that the odd numbers
are colored arbitrarily, and so that for each even number $2n$, the color
of $2n$ is different from the color of $n$. Such a coloring avoids monochromatic
$(a,2a)$-triples since such a triple has the form $\{x,y,z\}$ where $z=2y$.

We next consider the case $b>2a$.
Let $M=
4a(b^3+b^2-3b-3)+2b^3+4b^2+6b$ and let $\chi:[1,M] \rightarrow \{0,1\}$
be a $2$-coloring.
Assume there
is no monochromatic $(a,b)$-triple.  Then within the set
$\{2,4,\dots,2b+4\}$ there exist $x$ and $x+2$ that are not the
same color, since otherwise $\{2,2a+2,2b+4\}$ would be 
a monochromatic $(a,b)$-triple.  Without loss
of generality, assume $\chi(x)=0$ and 
$\chi(x+2)=1$.  Let $z$ be the least integer
greater than  $a(x+2)$ such that $b-2a$ divides
$z$.

  Let $S=\{z,az+(b-2a),bz+2(b-2a)\}$.  Since $z \leq 2a(b+1)+b$, we have
  \begin{equation}
   bz+2(b-2a) \leq 2a(b^{2}+b-2)+b^{2}+2b \leq M.
  \end{equation}
   Hence, since $S$
is an $(a,b)$-triple, some member, say $s$, of $S$ has color $1$.  Let 
\[ T = \{s+i(b-2a): 0 \leq i \leq \frac{s(b-1)}{b-2a} + 2 \} .\]
Note that $bs+2(b-2a)$ is the largest member of $T$, and that
$as + (b-2a) \in T$ since $b-2a$ divides $s$.
Also, by (1)
\begin{equation}
bs+2(b-2a) \leq 2a(b^{3}+b^{2}-2b-2) + b^{3} + 2b^{2} +2b \leq M,
\end{equation}
 so some member
of $T$ must have color 0 (otherwise $\{s,as+(b-2a),bs+2(b-2a)\}$ is
monochromatic).

  Let $t$ be the least member
of $T$ with color 0.  Then
$\chi(t-(b-2a))=1$. Note that (2) implies that
\[b(x+2)+2(t-ax-b) =2t+x(b-2a) \leq M.\]
 Thus, since
$\chi(x+2)=\chi(t-(b-2a))=1$, we must have
$\chi(b(x+2)+2(t-ax-b))=0$ (that $t - (b-2a) > a(x+2)$ follows from the
definition of $t$). This implies that
$\{x,t,bx+2(t-ax)\}$
is a monochromatic
$(a,b)$-triple, a contradiction.

\vspace{.1in}

The case for $b<2a$ is very similar.  Let
$M=
4a(b^3+2b^2+2b)-4b^2$ and let $\chi$ be a $2$-coloring of
$[1,M]$.  Then the set $\{2,4,\dots,2b+4\}$ contains 
$x-2$ and $x$ that are not the same color.
Assume
$\chi(x)=0$ and $\chi(x-2)=1$, and  let $z$ be
the least integer greater than $ax-(2a-b)$ such
that $2a-b$ divides $z$.  Let
$S=\{z,az+(2a-b),bz+2(2a-b)\}$. Let $s \in S$ have
color $1$ and define $T=\{s,s+(2a-b),s+2(2a-b),\dots,
bs+2(2a-b)\}$.  As in the previous case, $as+(2a-b) \in T$ and
$T \subseteq [1,M]$.  Hence, $T$ must have a least member, $t$, with
color $0$.  Then $\chi(t-(2a-b))=1$, and 
since $\chi(x-2)=\chi(t-(2a-b))=1$, we must have
$\chi(b(x-2)+2(t-ax+b))=0$.  This gives the monochromatic $(a,b)$-triple
$\{x,t,bx+2(t-ax)\}$, a contradiction. (That $bx+2(t-ax) \leq M$ follows easily
from the definitions of $z$, $s$, and $t$, and the fact that $x \geq 4$.)
\hfill $\Box$
\vskip 20pt

For certain pairs $(a,b)$, we are able to improve the upper bounds of Theorem
2.1. The next theorem deals with the case in which $a=b$. Theorem 2.1 gives
an upper bound for this case of $O(a^{4})$. The following theorem improves
this to $O(a^{2})$.

\vspace{.07in}

\noindent
{\bf Theorem 2.2:}
$
N(a,a;2) \leq
\left\{
\begin{array}{ll}
3a^2+a&\mathrm{for}\,\, 4 \leq a \,\, \mathrm{even}\\
8a^2+a&\mathrm{for}\,\, a \,\, \mathrm{odd}\\
\end{array}
\right.
$

\vspace{.07in}

\noindent
{\em Proof.}      
We start with the case when $a$ is even.  We may assume
that $a \geq 6$ since $N(4,4;2)=40$ was
obtained by computer search (for other 
exact values see Table 1 at the end of
this section).
We shall show that every red-blue coloring of $S=[1,3a^2+a]$
yields a monochromatic
$(a,a)$-triple by considering all possible 2-colorings of
the set $\{1,a+1,(3/2)a^{2}+a,2a^{2}+a\}$. Assume, by way of contradiction,
that there is a 2-coloring $\chi$ of $S$ that yields
no monochromatic $(a,a)$-triple.

Let $R$ be the set of red elements of $S$ under $\chi$,
and $B$ the set of blue elements of $S$ under $\chi$.
Without loss of generality we
assume $1 \in R$.

\vspace{.1in}

\noindent
{\tt Case I:} $a+1,2a^{2}+a \in R$.

\vspace{.1in}

\noindent
We then have the following implications.

\noindent $1,a+1 \in R \RA a+2 \in B$ (by considering the triple with $d=1$).\\
$1,2a^{2}+a \in R \RA a^{2}+a \in B$  (taking $d=a^{2}$).\\
$a+1,2a^{2}+a \in R \RA (3/2)a^{2}+a \in B$.\\
$a+2,(3/2)a^{2}+a \in B \RA (5/4)a^{2}+ (3/2)a \in R$.\\
$a^{2}+a,(3/2)a^{2}+a \in B \RA a/2+1 \in R$.\\
$a+1,(5/4)a^{2}+(3/2)a \in R \RA (3/2)a^{2}+2a \in B$.\\
$a/2+1,(5/4)a^{2}+(3/2)a \in R \RA 2a^{2}+2a \in B$.\\
This gives a contradiction since $\{a+2,(3/2)a^{2}+2a,2a^{2}+2a\}$ is a
monochromatic $(a,a)$-triple.

\vspace{.1in}

\noindent
{\tt Case II:} $a+1,(3/2)a^{2}+a \in R$ and $2a^{2}+a \in B$. 

\vspace{.1in}

\noindent
As in Case I, we must have $a+2 \in B$.
The following sequence of implications then leads to a contradiction. \\
$a+1,(3/2)a^{2}+a \in R \RA (5/4)a^{2}+a \in B$.\\
$a+2,2a^{2}+a \in B \RA 3a^{2} \in R$.\\
$1,(3/2)a^{2}+a \in R \RA 3a^{2}+a \in B$.\\
$a+2,3a^{2}+a \in B \RA 2a^{2}+(3/2)a \in R$.\\
$3a^{2},2a^{2}+(3/2) \in R \RA a+3 \in B$.\\
$a+2,(5/4)a^{2}+a \in B \RA (3/2)a^{2} \in R$.\\
$1,(3/2)a^{2} \in R \RA 3a^{2}-a \in B.$             \\
We now have a contradiction since $\{a+3,2a^{2}+a,3a^{2}-a\}$ is a blue
$(a,a)$-triple.

\vspace{.1in}

\noindent
{\tt Case III:} $(3/2)a^{2}+a,2a^{2}+a \in B$.

\vspace{.1in}

\noindent
This implies $a+1 \in R$, so that again we have $a+2 \in B$. Then \\
$a+2,2a^{2}+a \in B \RA (3/2)a^{2}+(3/2)a \in R$ and $3a^{2} \in R$. \\
$a+1,(3/2)a^{2}+(3/2)a \in R \RA 2a^{2}+2a \in B$ .\\
$a+2,(3/2)a^{2}+a \in B \RA 2a^{2} \in R$.\\
$2a^{2},3a^{2} \in R \RA a \in B$. \\
Hence, the $(a,a)$-triple $\{a,(3/2)a^{2}+a,2a^{2}+2a\}$ is blue, again
a contradiction.
\newpage

\noindent
{\tt Case IV:} $(3/2)a^{2}+a,2a^{2}+a \in R$. 

\vspace{.1in}

\noindent
Using this assumption and the fact that $1 \in R$, we have
$(3/4)a^{2}+a \in B$, $a^{2}+a \in B$, and $3a^{2}+a \in B$. Then \\
$(3/4)a^{2}+a,a^{2}+a \in B \RA (a/2)+1 \in R$. \\
$(a/2)+1,(3/2)a^{2}+a \in R \RA (5/2)a^{2}+a \in B$. \\
$3a^{2}+a,(5/2)a^{2}+a \in B \RA 2a+1 \in R$. 

\vspace{.1in}
\noindent
We now consider two subcases. 

\vspace{.1in}
\noindent
{\it Subcase} (i):  $2 \in R$.

\noindent
Then
$2,(3/2)a^{2}+a \in R \RA 3a^{2} \in B$. \\
$2,2a+1 \in R \RA 2a+2 \in B$.\\
Thus the $(a,a)$-triple
$\{2a+2,(5/2)a^{2}+a,3a^{2}\}$ is monochromatic. 

\vspace{.1in} 
\noindent 
{\it Subcase} (ii):  $2 \in B$

\noindent
$2,(3/4)a^{2}+a \in B \RA (3/2)a^{2} \in R$. \\
$2,a^{2}+a \in B \RA 2a^{2} \in R$. \\
$(3/2)a^{2},a/2+1 \in R \RA a^{2}+a/2 \in B$. \\
$(3/2)a^{2},2a^{2} \in R \RA a \in B$. \\
Thus the $(a,a)$-triple $\{a,(3/2)a^{2},2a^{2}\}$ is monochromatic.

\vspace{.1in}

\noindent
{\tt Case V:}  $a+1$, $2a^{2}+a \in B$.  

\vspace{.1in}

\noindent
In this case both $(3/2)a^{2}+a$ and $3a^{2}+a$ must be red, so that
the $(a,a)$-triple $\{1,(3/2)a^{2}+a,3a^{2}+a\}$ is monochromatic.

\vspace{.1in}

\noindent
{\tt Case VI:} $a+1$, $(3/2)a^{2}+a \in B$. 

\vspace{.1in}

\noindent
This assumption implies that $(5/4)a^{2}+a$ and $2a^{2}+a$ are red. Then \\
$1,2a^{2}+a \in R \RA a^{2}+a \in B$. \\
$a^{2}+a,(3/2)a^{2}+a \in B \RA a/2+1 \in R$. \\
Thus the $(a,a)$-triple $\{a/2+1,(5/4)a^{2}+a,2a^{2}+a\}$ is monochromatic.
\vskip 20pt
We now move onto the situation where $a$ is odd.
We may assume that $a \geq 5$ since $N(1,1;2)$ is
the van der Waerden number $w(3;2)$, which equals nine,
and $N(3,3;2)=39$ (see Table 1). 
Our method is very similar to that of the even case. Here
 we 2-color $T=[1,8a^2+a]$
and consider the various ways in which the set
$\{4a+1,5a^2+a,8a^2+a\}$ may be colored.
The following six cases cover all possibilities.

\vspace{.1in}

\noindent
{\tt Case I:} $4a+1,5a^{2}+a \in R$.

\vspace{.1in}

\noindent
In this case $6a^{2}+a \in B$ and, since $1 \in R$, $(5a+1)/2 \in B$.

\vspace{.1in}
\noindent
We consider two subcases.

\vspace{.1in}

\noindent
{\em Subcase} (i): $2 \in B$.

\noindent
Then 
$2,(5a+1)/2 \in B \RA 3a+1 \in R$.\\
$1,3a+1 \in R \RA 5a+2 \in B$.\\
$5a+2,6a^{2}+a \in B \RA 7a^{2} \in R$. \\
$1,7a^{2} \in R \RA (7a^{2}+a)/2 \in B$. \\
$(5a+1)/2,(7a^{2}+a)/2 \in B \RA (9a^{2}+a)/2 \in R$.\\
$(9a^{2}+a)/2, 5a^{2}+a \in R \RA 4a \in B$.\\
$2,4a \in B \RA 6a \in R$.\\
$6a,7a^{2} \in R \RA 8a^{2} \in B$.\\
$6a^{2}+a,8a^{2} \in B \RA 4a+2 \in R$. \\
$4a,8a^{2} \in B \RA 6a^{2} \in R$. \\
This gives the monochromatic $(a,a)$-triple $\{4a+2,5a^{2}+a,6a^{2}\}$.

\vspace{.1in}

\noindent
{\em Subcase} (ii): $2 \in R$.

\noindent
$2,5a^{2}+a \in R \RA (5a^{2}+3a)/2 \in B$.\\
$(5a+1)/2,(5a^{2}+3a)/2 \in B \RA (5a^{2}+5a)/2 \in R$.\\
$(5a^{2}+5a)/2,5a^{2}+a \in R \RA 4 \in B$.\\
$2,4a+1 \in R \RA 6a+2 \in B$.\\
$(5a+1)/2,(5a^{2}+3)/2 \in B \RA (5a^{2}+5a)/2 \in R$.\\
$1,(5a^{2}+5a)/2 \in R \RA 5a^{2}+4a \in B$.\\
$2,(5a^{2}+5a)/2 \in R \RA 5a^{2}+3a \in B$.\\
$5a^{2}+3a,5a^{2}+4a \in B \RA 5a+2 \in R$.\\
$4,6a+2 \in B \RA 8a+4 \in R$.\\
Thus, $\{2,5a+2,8a+4\}$ is a monochromatic $(a,a)$-triple.

\vspace{.1in}

\noindent
{\tt Case II:}
$4a+1,8a^{2}+a \in R$ and $5a^{2}+a \in B$.

\vspace{.1in}

\noindent
By using an obvious ``forcing'' argument (as in the previous cases)
on the following sequence of
$(a,a)$-triples, it is a routine exercise to show that the $(a,a)$-triple
 $\{1,3a,5a\}$ must be red: \hspace{.05in} $\{1,(5a+1)/2,4a+1\}$,
 $\{1,4a^{2}+a,8a^{2}+a\}$, $\{4a+1,6a^{2}+a,8a^{2}+a\}$,
 $\{(5a+1)/2,4a^{2}+a,(11a^{2}+3a)/2\}$, $\{4a+1,(11a^{2}+3a)/2,7a^{2}+2a\}$,
 $\{3a,5a^{2}+a,7a^{2}+2a\}$, $\{5a,6a^{2}+a,7a^{2}+2a\}$.

\vspace{.1in}

\noindent
{\tt Case III:}
$4a+1 \in R$ and $5a^{2}+a,8a^{2}+a \in B$.

\vspace{.1in}

\noindent
For this case we may use the $(a,a)$-triples $\{2a+1,5a^{2}+a,8a^{2}+a\}$,
 $\{1,2a+1,3a+2\}$, $\{3a+2,5a^{2}+a,7a^{2}\}$, $\{3a+2,(11a^{2}+3a)/2,8a^{2}
 +a\}$, $\{1,(7a^{2}+a)/2,7a^{2}\}$, $\{4a+1,(11a^{2}+3a)/2,7a^{2}+2a\}$,
  $\{2a,(7a^{2}+a)/2,5a^{2}+a\}$, $\{3a,5a^{2}+a,7a^{2}+2a\}$ to prove that
  the $(a,a)$-triple $\{1,2a,3a\}$ must be red.

  \vspace{.1in}

\noindent
{\tt Case IV:}
$5a^{2}+a,8a^{2}+a \in R$.

\vspace{.1in}

\noindent
By considering the triples $\{1,4a^{2}+a,8a^{2}+a\}$,
 $\{2a+1,5a^{2}+a,8a^{2}+a\}$, $\{2a+1,3a^{2}+a,4a^{2}+a\}$, and
 $\{2a+1,4a^{2}+a,6a^{2}+a\}$, we find that the $(a,a)$-triple
 $\{1,3a^{2}+a,6a^{2}+a\}$ must be red.

\vspace{.1in}

\noindent
{\tt Case V:}
$5a^{2}+a \in R$ and $4a+1,8a^{2}+a \in B$.

\vspace{.1in}

\noindent
In this case we have $6a^{2}+a \in R$ and hence $3a^{2}+a \in B$.

\vspace{.1in}
\noindent
We now consider two subcases.

\newpage
\noindent
{\em Subcase} (i): $2 \in B$.\\
By examining the triples $\{2,3a^{2}+a,6a^{2}\}$, $\{4a+2,5a^{2}+a,6a^{2}\}$,
$\{6a-1,6a^{2},6a^{2}+a\}$, $\{6a-1,7a^{2},8a^{2}+a\}$,  $\{2,3a+1,4a+2\}$,
$\{3a+1,4a^{2}+a,5a^{2}+a\}$, $\{2,4a^{2}+a,8a^{2}\}$, $\{4a,6a^{2},8a^{2}\}$,
$\{6a,7a^{2},8a^{2}\}$, we find that the $(a,a)$-triple $\{2,4a,6a\}$ must be blue.

\vspace{.1in}
\noindent
{\em Subcase} (ii): $2 \in R$. \\
By examining the triples $\{2,(5a^{2}+3a)/2,5a^{2}+a\}$,
$\{2a+2,(5a^{2}+3a)/2,3a^{2}+a\}$, $\{2,2a+2,2a+4\}$, $\{2,2a+1,2a+2\}$,
$\{2a+2,(7a^{2}+3a)/2,5a^{2}+a\}$, $\{2a+4,(5a^{2}+5a)/2,3a^{2}+a\}$,
$\{(3a+3)/2,(5a^{2}+3a)/2,(7a^{2}+3a)/2\}$, $\{2a+1,(5a^{2}+3a)/2,3a^{2}+2a\}$,
$\{1,(5a^{2}+5a)/2,5a^{2}+4a\}$, $\{(3a+3)/2,3a^{2}+2a,(9a^{2}+5a)/2\}$,
 we find that the $(a,a)$-triple $\{4a+1,(9a^{2}+5a)/2,5a^{2}+4a\}$ is blue.

\vspace{.1in}
\noindent
{\tt Case VI:}
$4a+1,5a^{2}+a \in B$.

\vspace{.1in}

\noindent
The sequence of triples $\{4a+1,5a^{2}+a,6a^{2}+a\}$, $\{1,3a^{2}+a,6a^{2}+a\}$,
 $\{a+1,3a^{2}+a,5a^{2}+a\}$, $\{1,a+1,a+2\}$, $\{a+2,3a^{2}+a,5a^{2}\}$,
 $\{4a-1,5a^{2},6a^{2}+a\}$, $\{1,(5a^{2}+a)/2,5a^{2}\}$, 
 $\{2a,(5a^{2}+a)/2,3a^{2}+a\}$, $\{a+2,(5a^{2}+a)/2,4a^{2}-a\}$, 
 $\{2a,4a^{2}-a,6a^{2}-2a\}$, $\{1,4a^{2}-a,8a^{2}-3a\}$ leads us to
 conclude that the $(a,a)$-triple $\{4a-1,6a^{2}-2a,8a^{2}-3a\}$ is blue.
\hfill $\Box$
\vskip 20pt

Another circumstance for which we can improve the upper bounds of Theorem 2.1
is the case in which $b=2a-1$ (for $a=1$ this is the van der Waerden
number $w(3,2)$). By Theorem 2.1, $N(a,2a-1;2)$ is bounded above by a
function having order of magnitude $32a^{4}$. We can improve this
to $16a^3$ by making
use of the following theorem which is taken from  [1].

First, we introduce some notation.
Let $f: {\bf N} \rightarrow {\bf R^{+}}$ be a non-decreasing function. Denote
by $w(f,k)$ the least positive integer (if it exists) such that whenever
$[1,w(f,k)]$ is 2-colored, there must exist a monochromatic $k$-term
arithmetic progression $\{a,a+d,a+2d,\ldots,a+(k-1)d\}$ with $d \geq f(a)$.
In [1] it is shown that $w(f,3)$ always exists, and bounds for this function
are given as follows.

\vspace{.1in}

\noindent
{\bf Theorem 2.3:} (Brown and Landman [1])
Let $f:\mathbf{N} \rightarrow \mathbf{R^+}$ be a non-decreasing
function.  Let $b=1+ 4 \lceil \frac{f(1)}{2} \rceil$.  Then
\[w(f,3) \leq \lceil 4f(b+4 \lceil \frac{f(b)}{2} \rceil)+14   
\lceil \frac{f(b)}{2} \rceil + 7b/2 - 13/2 \rceil.\]  Further,
if $f$ maps into $\mathbf{N}$ with $f(n) \geq n$ for all
$n \in \mathbf{N}$, then
$ w(f,3,2) \geq 8f(h)+2h+2-c$, 
where $h=2f(1)+1$ and $c$ is the largest integer such that
$f(c)+c \leq 4f(h)+h+1$.

\vspace{.15in}

Relating Theorem 2.3 to $(a,b)$-triples, we have the following corollary.

\vspace{.1in}

\noindent
{\bf Corollary 2.1:}  
For all $a \geq 2$,
$$
16a^2-12a+6 \leq
N(a,2a-1;2) \leq   
\left\{
\begin{array}{ll}
16a^3-2a^2+4a-3&\mathrm{for} \,\, a \,\, \mathrm{even}\\
16a^3+14a^2+2a-3&\mathrm{for} \,\, a \,\, \mathrm{odd}\\
\end{array}
\right.
$$
\newpage
\noindent
{\em Proof.}  Note that $\{x,y,z\}$ is an $(a,2a-1)$-triple if and
only if it is an arithmetic progression with $y-x \geq (a-1)x+1$.
By applying Theorem 2.3 
with $f(x)=(a-1)x+1$ 
we obtain the desired bounds.
\hfill $\Box$
\vskip 20pt
We now present some lower bounds for all
$(a,b)$-triples.  This is done by providing
$2$-colorings which avoid monochromatic
$(a,b)$-triples.

\vspace{.1in}

\noindent
{\bf Theorem 2.4:} If $b \geq 2a$ then 
$N(a,b;2) \geq 2b^2+5b-(2a-4)$.
If $b < 2a$ then
$N(a,b;2) \geq 
3b^2-(2a-5)b - (2a-4)$.

\vspace{.07in}

\noindent
{\em Proof.} For the case $b \geq 2a$, we  
will exhibit a $2$-coloring of $[1,2b^2+5b-2a+3]$ with no
monochromatic $(a,b)$-triple.  Color
$[b+2,b^2+2b+1]$ red and its complement blue.  It is an easy
exercise to show that monochromatic $(a,b)$-triples are avoided.

For the case where $b<2a$ the $2$-coloring of
$[1,3b^2-(2a-5)b-(2a-4)]$ with $[b+2,b^2+2b+1]$ colored red
and its complement colored blue is easily seen to avoid
monochromatic $(a,b)$-triples. 
\hfill $\Box$
\vskip 20pt

We are able to improve slightly the lower
bound given in Theorem 2.5
for the case when $a=1$.
In fact, from computer calculations
(see Table 1 below), it appears that this
inequality may in fact be an equality.

\vspace{.1in}

\noindent
{\bf Theorem 2.5:} $N(1,b;2) \geq 2b^2+5b+6$ for all $b \geq 3$.

\vspace{.07in}

\noindent
{\em Proof.}
Consider the following red-blue coloring of $[1,2b^2+5b+5]$:
color $[1,b+1]$, $\{b+3\}$, and $[b^2+2b+4,2b^2+5b+5]$
red and the other integers blue.  We now show that
this coloring avoids monochromatic $(1,b)$-triples.

 Assume $\{x,y,z\}=\{x,x+d,bx+2d\}$ is a blue $(1,b)$-triple.
 Since the largest blue element is
$b^2+2b+3$, we must have $x=b+2$.  Thus, since we
must have $d \geq 2$, we see that 
$z > b^2+2b+3$, which is not possible.

Now assume $\{x,y,z\}$ is red.
First, if $x,x+d \in \{1,2,\dots,b,b+1,b+3\}$ then
we have $b+2 \leq bx+2d \leq b^2+b+4$.  Hence, 
the only possibility here is $z=b+3$, but
$bx+2d=b+3$ has no solution in $x$ for $b \geq 3$.
Second, if $y \in [b^2+2b+4,2b^2+5b+5]$ then
$bx+2d \geq bx + 2b^2+2b+2$.  Hence, we must
have $x \in \{1,2,3,4\}$ ($4$ is possible if $b=3$).
However, this gives $bx+2d \geq bx + 2b^2+4b$, 
which implies that $x \in \{1,2\}$ ($2$ is possible
if $b=3$).  This in turn implies that
$z \geq bx+2b^2+2b+4$ which gives $x=1$ as
the only possibility.  However with $x=1$ we must
have $z > 2b^2+5b+5$, which is out of bounds.
\hfill $\Box$
\vskip 20pt

Below we present a table of computer-generated values for $N(a,b;2)$ for
small $a$ and $b$. We also include  computer-generated lower bounds
 for those cases
where the computer time became excessive
 (the program is available for download as the
Fortran77 program {\em VDW.f} at
{\tt http://math.colgate.edu/\~{}aaron/}).

\newpage
\begin{center}
$\mathbf{N(a,b;2)}$ {\bf Values}
\end{center}

\begin{center}
\begin{tabular}{|l||c|c|c|c|c|c|c|} \hline
$a \setminus b$&$1$&$2$&$3$&$4$&$5$&$6$&$7$\\ \hline
$1$&9&dne&$39$&$58$&$81$&$\geq 108$&$\geq 139$\\ \hline
$2$&&$16$&$46$&dne&$139$&$\geq 106$&$\geq 133$\\ \hline
$3$&&&$39$&$60$&$114$&dne&$\geq 135$\\ \hline
$4$&&&&$40$&$87$&$\geq 124$&$\geq 214$\\ \hline
$5$&&&&&$70$&$100$&$\geq 150$\\ \hline
$6$&&&&&&$78$&$\geq 105$\\ \hline
$7$&&&&&&&$95$\\ \hline
\end{tabular}
\end{center}
\begin{center}
{\bf Table 1}
\end{center}
\vskip50pt

\begin{center}
{\bf 3. The Degree of Regularity of $\mathbf{(a,b)}$}
\end{center}

In this section we consider $N(a,b;r)$ for general $r$. We begin by
showing, in Theorems 3.1 and 3.2, that for many choices of $a$ and $b$, the
pair $(a,b)$ is not regular.
For such pairs we find an upper bound on $dor (a,b)$.

\vspace{.1in}

\noindent
{\bf Theorem 3.1:} Let $1 \leq a < b$, and assume that
$b \geq (2^{3/2}-1)a + 2 - 2^{3/2}$.  Let 
$c = \lceil b/a \rceil$.
 Then $dor(a,b) \leq  \lceil log_{\sqrt{2}}\,\, c \rceil$.

\vspace{.07in}

\noindent
{\em Proof.}
Let $r = \lceil \log_{\sqrt{2}} c \rceil + 1$.
We will give an $r$-coloring of the natural numbers which
contains no monochromatic $(a,b)$-triple. 
For readability, let $p = \sqrt{2}$.
 Using the colors
$0,1, \dots, r-1$, define the coloring $\chi$ by letting
$\chi(x) \equiv i$ (mod $r$), where $p^i \leq x < p^{i+1}$.

Assume that there exists an $(a,b)$-triple, say $x<y<z$, that is
monochromatic under $\chi$.  Let $j$ be the integer such that
$p^j \leq y < p^{j+1}$.  Since $\{x,y,z\}$ is an $(a,b)$-triple, 
$y=ax+d$ and $z=bx+2d$ for some $d$.  Thus
$ z \leq cy < p^{r-1} p^{j+1} = p^{j+r}$.  Hence, by the way
$\chi$ is defined and the fact that $\chi(y)=\chi(z)$, we must have 
$p^{j} \leq y < z < p^{j+1}$.

We now consider two cases. 

\vspace{.1in}
\noindent
{\tt Case I:} $b \geq 2a-1$.

In this case, $y-x = (a-1)x + d \leq (b-a)x + d = z-y < 
p^{j}(p-1) \leq p^j(1 - 1/p^{r-1})$.
Hence, since $y > p^j$, we have that $x \geq p^j - p^j(1-1/p^{r-1}) 
= p^{j-r+1}$.  Since $\chi(x)=\chi(y)$, and by the definition of $\chi$, 
we must have $p^j \leq x<y<p^{j+1}$.  Thus we have that all three
numbers $x,y,z$ belong to the interval $[p^j,p^{j+1})$.  Hence,
$z-x = (b-1)x + 2d < p^j(p-1) \leq x(p-1)$, a contradiction
(since $b-1 > p-1$).
\vspace{.1in}

\noindent
{\tt Case II:}  $c=2$ and
$b \geq (2^{3/2}-1)a + 2 - 2^{3/2}$.

In this case, $2(a-1) \leq (b-a)/(p-1)$, so that
$(a-1)x/(b-a)x \leq 1/(2p-2)$.  Therefore,
$((a-1)x+d)/((b-a)x+d) \leq 1/(2p-2)$.  Hence,
$y-x \leq (z-y)/(2p-2) < (p-1)p^j/(2p-2) = p^j/2 = p^{j-2}$.
So, $x \geq p^j - p^{j-2} = p^{j-2}$.
Since $r=3$ in this case, and
$\chi(x)=\chi(y)$, we must have $p^j \leq x < p^{j+1}$.
Thus, as in Case I, $x$ and $z$ both belong to
the interval $[p^j,p^{j+1})$, and we again have 
a contradiction.
\hfill $\Box$
\vskip 20pt

In the following theorem we give an upper bound on $dor(a,b)$ for several
pairs $(a,b)$ that are either not covered by Theorem 3.1 or for which we
are able to improve the bound of Theorem 3.1.

\vspace{.1in}
 
\noindent
{\bf Theorem 3.2:} $dor(1,3) \leq 3$, $dor(2,2) \leq 5$, $dor(2,5) \leq 3$, 
$dor(2,6) \leq 3$, $dor(3,3) \leq 5$, $dor(3,4) \leq 5$, 
$dor(3,8) \leq 3$, and $dor(3,9) \leq 3$. 

\vspace{.07in}

\noindent
{\em Proof.} We give the proof for the pair $(2,2)$, and outline the
proofs for the other cases, which are quite similar.

To show that $dor(2,2) \leq 5$, we provide a 6-coloring of the positive
integers that avoids monochromatic $(2,2)$-triples.
Let
\[ \chi(i)= \left\{
\begin{array}{ll}
2k \,\, (\mathrm{mod}\,\, 6) & \mathrm{if}\,\, i \in [2^k,\lfloor p2^k \rfloor)\\
2k+1 \,\, (\mathrm{mod}\,\, 6) & \mathrm{if}\,\, i \in [\lfloor p2^k \rfloor,2^{k+1})
\end{array}
\right.\]
where $p = \sqrt{2}$ .
Assume $\{x,2x+d,2x+2d\}$ is a (2,2)-triple such that $\chi(x)=\chi(2x+d)$.
We will show that $\chi(2x+2d) \neq \chi(x)$. We consider two cases.

\vspace{.1in}
\noindent
{\tt Case I:} $2^k \leq x < \lfloor p2^k \rfloor$
              for some $k \in \{0,1,2, \dots\}$.

\noindent
Since $\chi(x)=\chi(2x+d)$ and $2x+d > p2^{k}$, there exists an $m \in {\bf N}$
such that $2x+d \in [2^{k+3m}, \lfloor p2^{k+3m} \rfloor)$.
Hence $d \in [2^{k+3m}-p2^{k+1}, p2^{k+3m}-2^{k+1}]$. This yields
\begin{equation}
 2x+2d \geq 2^{k+3m}+2^{k+3m}-p2^{k+1} \geq p2^{k+3m},
 \end{equation}
 and
 \begin{equation}
 2x+2d < p2^{k+3m+1}-2^{k+1} < 2^{k+3(m+1)}.
 \end{equation}
 By (3) and (4) it follows that $\chi(2x+2d) \neq \chi(x)$.

 \vspace{.1in}

 \noindent
 {\tt Case II:}  
$\lfloor p2^k \rfloor \leq x < 2^{k+1}$
               for some $k \in \{0,1,2, \dots \}$.

\noindent

As in Case I, there must exist an $m \in \{1,2,\dots\}$ such that
$2x+d \in [\lfloor p2^{k+3m} \rfloor, 2^{k+3m+1})$.
Thus,
$d \in [p2^{k+3m}-2^{k+2}-1,
2^{k+3m+1}-p2^{k+1}]$.    Therefore
\begin{equation}
2x+2d \geq 2^{k+3m+1}(p-2^{1-3m}-2^{-k-3m}) \geq 2^{k+3m+1},
\end{equation}
and
\begin{equation}
2x+2d < 2^{k+3m+2}-p2^{k+1} < 2^{k+3(m+1)}.
\end{equation}
It follows from (5) and (6) that $\chi(2x+2d) \neq \chi(x)$.

\vspace{.07in}

The proofs that $dor(3,3) \leq 5$ and $dor(3,4) \leq 5$  may be done in the
same way as that for $dor(2,2)$ except that we use $p=\sqrt{3}$ instead
of $p=\sqrt{2}$ and we use powers of $3$ instead of $2$ in the
defined intervals.

The cases of $(2,5)$, $(2,6)$, $(3,8)$, and $(3,9)$ are done
similarly, where we use
a 4-coloring rather than a 6-coloring, which is defined the same as $\chi$
except that ``mod 6'' is replaced by ``mod 4;'' where we take $p$ to be
$1.6$, $1.5$, $1.9$, and $1.9$, respectively; and where the powers in the
defined intervals are powers of the given value of $a$.  The case
$(1,3)$ is done using a ``mod 4" coloring with $p=\sqrt{3}$,
where the
powers in the given intervals of the coloring are powers of $3$.
\hfill $\Box$
\vskip 20pt

By using Theorems 2.4 and 2.5, we are able to obtain the following lower
bounds for $N(a,b;r)$.

\vspace{.1in}

\noindent
{\bf Proposition 3.1:}  If $b \geq 2a$ then
$N(a,b;r) \geq 2b^r+5b^{r-1}-(2a-4)b^{r-2}+
\sum _{i=0}^{r-3} \,\, b^i$ for $r \geq 2$.
If $b<2a$ then
$N(a,b;r) \geq 3b^r-(2a-5)b^{r-1}-(2a-4)b^{r-2}+2   
\sum_{i=0}^{r-3}\,\,b^i$ for $r \geq 2$.

\vspace{.07in}

\noindent
{\em Proof.} We induct on $r$.  The case
$b \geq 2a$ and $r=2$ is proved in Theorem 2.4. 
Assuming $r \geq 3$ and that the result holds
for $r-1$ with $b \geq 2a$, there exists an $(r-1)$-coloring
of $[1,2b^{r-1}+5b^{r-2}-(2a-4)b^{r-3}+2\sum_{i=0}^{r-4} \,\, b^i
\,\, -1]$ with no monochromatic triple.  Color the
interval  
$[2b^{r-1}+5b^{r-2}-(2a-4)b^{r-3}+2\sum_{i=0}^{r-4} \,\, b^i,
2b^{r}+5b^{r-1}-(2a-4)b^{r-2}+2\sum_{i=0}^{r-3} \,\, b^i \,\, -1]$
with the remaining color.  By construction this
$r$-coloring avoids monochromatic triples.

The case $b<2a$ is quite similar and will be omitted.
\hfill $\Box$
\vskip 20pt

\noindent
{\bf Proposition 3.2:} $N(1,b;r) \geq 
2b^r+5b^{r-1}+6b^{r-2}+2 \sum_{i=0}^{r-3}\,\,b^i$ for all
$b,r \geq 2$.

\vspace{.07in}

\noindent
{\em Proof.} We induct on $r$.  The case
$r=2$ is proved in Theorem 2.5.  Assuming $r \geq 3$ and that the
inequality is true for $r-1$, we have the existence
of an $(r-1)$-coloring of
$[1,2b^{r-1}+5b^{r-2}+6b^{r-3}+2 \sum_{i=0}^{r-4}\,\,b^i\,\,-1]$
which does not contain a monochromatic
$(1,b)$-triple.  Color the interval
$[2b^{r-1}+5b^{r-2}+6b^{r-3}+2 \sum_{i=0}^{r-4}\,\, b^i,
2b^r+5b^{r-1}+6b^{r-2}+2 \sum_{i=0}^{r-3}\,\,b^i-1]$ 
with the remaining color.  It is an easy exercise to show that
there is no monochromatic $(1,b)$-triple in this $r$-coloring. 
\hfill $\Box$

\vspace{.2in}

We conclude this section with a table which describes what is known about
$dor(a,b)$ for some small values of $a$ and $b$. By Theorem 2.1, we know that
if $b \neq 2a$, then $dor(a,b) \geq 2$. In the third column of the table
below we give the reason for the given upper bound on $dor(a,b)$.

\newpage
\begin{center}
{\bf Values of $\mathbf{dor(a,b)}$}
\end{center}

\begin{center}
\begin{tabular}{|c|c|l|} \hline
$(a,b)$ & $dor(a,b)$ &               reason     \\ \hline
(1,1) & $\infty$     & van der Waerden's Theorem \\
(1,2) & $1$          & Theorem 2.1                 \\
(1,3) & $2-3$          & Theorem 3.2            \\
(1,4) & $2-4$         & Theorem 3.1               \\
(1,5) & $2-5$         & Theorem 3.1 \\
(1,6) & $2-6$  &      Theorem 3.1 \\
(1,7) & $2-6$ &       Theorem 3.1 \\
(1,8) & $2-6$  &      Theorem 3.1 \\
(1,9) & $2-7$   &     Theorem 3.1 \\
(2,2) & $2-5$ &       Theorem 3.2\\
(2,3) & $2$         &   Theorem 3.1\\
(2,4) & $1$       &     Theorem 2.1 \\
(2,5) & $2-3$     &     Theorem 3.2 \\
(2,6) & $2-3$      &    Theorem 3.2 \\
(2,7) & $2-4$      &    Theorem 3.1 \\
(2,8) & $2-4$     &     Theorem 3.1 \\
(2,9) & $2-5$    &       Theorem 3.1\\
(3,3) & $2-5$ &        Theorem 3.2 \\
(3,4) & $2-5$ &         Theorem 3.2 \\
(3,5) & $2$   &         Theorem 3.1 \\
(3,6) & $1$   &         Theorem 2.1 \\
(3,7) & $2-4$  &         Theorem 3.1 \\
(3,8) & $2-3$  &         Theorem 3.2 \\
(3,9) & $2-3$  &         Theorem 3.2 \\ \hline
\end{tabular}
\end{center}
\begin{center}
{\bf Table 2}
\end{center}

\vspace{.2in}

\begin{center}
{\bf 4. A More General Question}
\end{center}

In this section we move from $(a,b)$-triples to
sets of the form
$\{x,ax+d,ax+2d,\dots,ax+(k-1)d\}$ for $a \geq 1$ and
$k \geq 3$.  Let us call such a set a {\em $k$-term $a$-progression}.
For $a=1$ these are simply the $k$-term
arithmetic progressions.  Van der Waerden's theorem
states that given $r \geq 1$, any $r$-coloring of
the natural numbers must contain arbitrarily
long monochromatic arithmetic progressions.
Theorem 4.1 shows that a similar result
does {\it not} hold for $a>1$ and $r>3$.
 Denote by $dor_{k}(a)$ the largest number
of colors with which we can arbitrarily color {\bf N} and be guaranteed
the existence of a monochromatic $k$-term $a$-progression. Theorem 4.1
shows that for $k$ large enough, $dor_{k}(a) \leq 3$. 
\newpage
\noindent
{\bf Theorem 4.1:} For all $a \geq 2$ and all integers 
$k \geq \frac{a^2}{a+1}+2$, $dor_{k}(a) \leq 3$.

\vspace{.07in}

\noindent
{\em Proof.}  It suffices to exhibit a $4$-coloring
of $\mathbf{N}$ which avoids monochromatic $k$-term $a$-progressions.
Clearly, we may assume $k= \lceil
\frac{a^2}{a+1}\rceil +2$.

Define a $4$-coloring of $\mathbf{N}$ by
coloring each interval $[a^j,a^{j+1})$ with the color $j$ (mod $4$). 
We will show that there is no monochromatic $k$-term $a$-progression
by showing that if $x$ and $ax+d$ are the same color,
then $ax+(k-1)d$ is a different color.

Let $x \in [a^i,a^{i+1})$, and assume $ax+d$ has the same color as $x$.
  Then clearly
$ax+d \not \in [a^i,a^{i+1})$.  Hence, there exists
an $m \in \mathbf{N}$ such that
$ax+d \in [a^{i+4m},a^{i+4m+1})$.  From this we
conclude that
\[a^{i}(a^{4m}-a^2) \leq d \leq a^{i+1}(a^{4m}-1).\]
To complete the proof we will show that 
 \begin{equation}
 ax+(k-1)d < a^{i+4(m+1)}
\end{equation}
and
\begin{equation}
a^{i+4m+1} \leq ax+(k-1)d.
\end{equation}
From (7) and (8) we can conclude that $ax+(k-1)d$ is colored differently than
$x$ and $ax+d$.

To prove (7), note that $k < a^3+1$ for all $a \geq 2$.
Thus
\[ 1+(k-2)(1-a^{-4m})<a^3, \] and hence
\[ a^{i+4m+1} +(k-2)a^{i+1}(a^{4m}-1) <  a^{i+4(m+1)}.\]
This last inequality, together with the fact that
\[ ax+(k-1)d = ax+d+(k-2)d \leq a^{i+4m+1}+(k-2)a^{i+1}(a^{4m}-1),\]
implies (7).

To prove (8), first note that since $k \geq a^2/(a+1)+2$, we have
$(k-2)(a^2-1) \geq a^3-a^2$. Hence
\[ a^{i+4m} + (k-2)a^{i}(a^{4m}-a^2) \geq a^{i+4m+1}.\]
This last inequality, along with the fact that
\[ ax+(k-1)d \geq a^{i+4m} + (k-2)a^{i}(a^{4m}-a^2), \]
shows that (8) holds.
\hfill $\Box$
\vskip 20pt

According to Theorem 4.1, it is not true that every 4-coloring of ${\bf N}$
yields arbitrarily long monochromatic $a$-progressions. We are not sure
 if this holds for two or three colors.  However, if
 for $r=2$ or $r=3$, every $r$-coloring of ${\bf N}$ does 
yield arbitrarily long monochromatic $a$-progressions, then a
somewhat
stronger result holds, as stated in Proposition 4.1 below.
We omit the proof, a trivial generalization of the proof of
[3, Theorem 2, p. 70].

\vspace{.1in}

\noindent
{\bf Proposition 4.1} Let $a \in {\bf N}$ and let $r \in \{2,3\}$. If for
every $r$-coloring of ${\bf N}$ there are arbitrarily long monochromatic
$a$-progressions, then for all $s \geq 1$ there exists $n=n(a,r,s)$ such that
if $[1,n]$ is $r$-colored then for all $k \in {\bf N}$ there exists
 $\hat{x},\hat{d}$
so that
 $\{\hat{x},a\hat{x}+\hat{d},a\hat{x}+2\hat{d},\ldots,a\hat{x}+k\hat{d}\}
 \cup \{s\hat{d}\}$
 is monochromatic.
\vskip 40pt
\begin{center}
{\bf 5. Some Concluding Remarks}
\end{center}

Although we have not proved that $dor(a,b)< \infty$ for general
$a$ and $b$, the evidence in this paper leads us to
believe that this is the case for all $(a,b) \neq (1,1)$.  In particular,
we make the following conjecture:

\vspace{.1in}

\noindent
{\bf Conjecture:} Let $a > 1$ and $r > 3$. Define $K(a,r)$ to be the
least positive integer such that $dor_{K}(a) \leq r$. Then there exists
an $s > r$ such that $K(a,s) < K(a,r)$.

\vspace{.1in}

By Theorem 4.1, we know that $K(a,r)$ exists.
Clearly, $K(a,s) \leq K(a,r)$ for $s \geq r$,
but if we are able to show that the inequality is strict for some $s$,
 then we can conclude that $dor(a,a) < \infty$ for all $a > 1$.
In fact it may be true that $dor(a,b)=2$ for all $b \neq 2a$, although we have
presented scant evidence for this.
\vskip 50pt
\noindent
{\bf References}

\noindent
[1] T.C. Brown and B. Landman, Monochromatic arithmetic
progressions with large differences, {\em Bull. Australian Math. Soc.}
{\bf 60} (1999), 21-35.

\noindent
[2] T.C. Brown, B. Landman, M. Mishna, Monochromatic homothetic copies of
$\{1,1+s,1+s+t\}$, {\em Canadian Math. Bull.} {\bf 40} (1997), 149-157.

\noindent
[3] R.L. Graham, B.L. Rothschild, J.H. Spencer, {\em Ramsey Theory, Second
Ed.},
John Wiley and Sons, New York, 1990.

\noindent
[4] B. Landman, Avoiding arithmetic progressions (mod $m$) and arithmetic
progressions, {\em Utilitas Math.} {\bf 52} (1997), 173-182.

\noindent
[5] B. Landman, Ramsey functions for quasi-progressions, {\em Graphs and
Combinatorics} {\bf 14} (1998), 131-142.

\noindent
[6] B. L. van der Waerden, Beweis einer Baudetschen
Vermutung, {\em Nieuw Arch. Wisk.} {\bf 15} (1927), 212-216.
\end{document}